\documentclass[final,1p,times]{elsarticle}

\usepackage{amssymb}
 \usepackage{amsthm}
\usepackage{amscd}
\usepackage{amsmath}
\usepackage{amsfonts}
\usepackage{amssymb}
\usepackage{graphicx}
\usepackage{color}
\newtheorem{theorem}{Theorem}
\newtheorem{proposition}[theorem]{Proposition}
\newtheorem{lemma}[theorem]{Lemma}

\usepackage{mathrsfs}
\usepackage{titletoc}
\newcommand{\D}{\Delta}
\newcommand{\ra}{\rightarrow}
\newcommand{\p}{\partial}
\newcommand{\f}{\frac}

\newcommand{\be}{\begin{equation}}
\renewcommand{\ra}{\rightarrow}
\newcommand{\ee}{\end{equation}}
\newcommand{\bea}{\begin{eqnarray}}
\newcommand{\eea}{\end{eqnarray}}
\newcommand{\bna}{\begin{eqnarray*}}
\newcommand{\ena}{\end{eqnarray*}}

\renewcommand{\le}{\left}
\newcommand{\ri}{\right}

\journal{***}

\begin{document}

\begin{frontmatter}

\title{Mean field equations on a closed Riemannian surface with the action of an
isometric group}

\author{Yunyan Yang}
 \ead{yunyanyang@ruc.edu.cn}
\author{Xiaobao Zhu}
 \ead{zhuxiaobao@ruc.edu.cn}

 \address{Department of Mathematics,
Renmin University of China, Beijing 100872, P. R. China}

\begin{abstract}
  Let $(\Sigma,g)$ be a closed Riemannian surface, $\textbf{G}=\{\sigma_1,\cdots,\sigma_N\}$ be an isometric group
  acting on it. Denote a positive integer $\ell=\inf_{x\in\Sigma}I(x)$, where $I(x)$ is the number of all distinct points of
  the set $\{\sigma_1(x),\cdots,\sigma_N(x)\}$. A sufficient condition for existence of solutions to the mean field equation
  $$\Delta_g u=8\pi\ell\le(\f{he^u}{\int_\Sigma he^udv_g}-\f{1}{{\rm Vol}_g(\Sigma)}\ri)$$
  is given. This recovers results of Ding-Jost-Li-Wang (Asian J Math 1997) when $\ell=1$ or equivalently
  $\textbf{G}=\{Id\}$, where $Id$ is the identity map.
\end{abstract}


\end{frontmatter}

%\tableofcontents
\section{Introduction}
Let $(\mathbb{S}^2,g_0)$ be the 2-dimensional sphere $x_1^2+x_2^2+x_3^2=1$ with the metric $g_0=dx_1^2+dx_2^2+dx_3^2$
and the corresponding volume element $dv_{g_0}$. It was proved by Moser \cite{Moser}, via a symmetrization argument, that
there exists some constant $C$ such that if
$u$ satisfies $\int_{\mathbb{S}^2}|\nabla_{g_0} u|^2dv_{g_0}\leq 1$ and $\int_{\mathbb{S}^2}udv_{g_0}=0$, then
\be\label{Mosr-2}\int_{\mathbb{S}^2} e^{\alpha u^2}dv_{g_0}\leq C,\quad\forall\alpha\leq 4\pi;\ee
moreover, $4\pi$ is the best constant in the sense that if $\alpha>4\pi$, then the integerals in (\ref{Mosr-2}) are still finite,
but there is no uniform constant $C$ such that (\ref{Mosr-2}) holds.
Also, it was pointed out by Moser \cite{Moser-73} that the same argument in \cite{Moser} indicates a similar inequality for even functions.
 Namely, there exists a constant $C$ such that if $u$ satisfies $u(-x)=u(x)$ for all $x\in \mathbb{S}^2$, $\int_{\mathbb{S}^2}|\nabla_{g_0}u|^2dv_{g_0}\leq 1$,
and $\int_{\mathbb{S}^2}udv_{g_0}=0$, then
\be\label{Mosr-even}
\int_{\mathbb{S}^2} e^{\alpha u^2}dv_{g_0}\leq C,\quad\forall\alpha\leq 8\pi.\ee
Observing that the inequality (\ref{Mosr-even}) can be applied to the prescribing Gaussian curvature equation, Moser \cite{Moser-73}
obtained the following: If $K:\mathbb{S}^2\ra \mathbb{R}$ is a smooth function satisfying $K(-x)=K(x)$ for all $x\in\mathbb{S}^2$, then there is a
smooth solution to
\be\label{pres-S2}\Delta_{g_0}u+1-Ke^{2u}=0,\ee
where $\Delta_{g_0}$ is the Laplace-Beltrami operator. The geometric meaning of the above equation is as follows: If $u$ is a smooth solution to
(\ref{pres-S2}), then the metric $g=e^{2u}g_0$ has the Gaussian curvature $K$.

Moser's inequality (\ref{Mosr-2}) was generalized by Fontana \cite{Fontana} to Riemannian manifolds case. In particular, let $(\Sigma,g)$ be a
closed Riemannian surface. Then for all smooth functions $u$ with $\int_\Sigma|\nabla_gu|^2dv_g\leq 1$ and $\int_\Sigma udv_g=0$, there exists
a constant depending only on $(\Sigma,g)$ such that
$$\int_\Sigma e^{4\pi u^2}dv_g\leq C;$$
moreover, $4\pi$ is the best constant for the above inequality.
Also Moser's inequality (\ref{Mosr-even}) was generalized to various versions. Let $\textbf{G}=
\{\sigma_1,\cdots,\sigma_N\}$ be an isometric group acting on a closed Riemannian surface $(\Sigma,g)$, $I(x)$ be the number of all distinct elements in $\textbf{G}(x)=\{\sigma_1(x),\cdots,\sigma_N(x)\}$, and $\ell=\inf_{x\in\Sigma}I(x)$. Using isoperimetric inequalities, among other results Chen \cite{Chen-90} proved that there exists a constant $C$ such that for all smooth functions $u$ with
$u(\sigma(x))=u(x)$ for all $\sigma\in\textbf{G}$ and all $x\in\Sigma$, $\int_\Sigma|\nabla_gu|^2dv_g\leq 1$ and $\int_\Sigma udv_g=0$, there holds
\be\label{Chen}\int_\Sigma e^{4\pi\ell u^2}dv_g\leq C.\ee
Recently, motivated by \cite{Yang-Tran, Lu-Yang-1, Yang-JDE-15}, Fang-Yang \cite{Fang-Yang} employed blow-up analysis to improve (\ref{Chen}) to
analogous inequalities involving eigenvalues of the Laplace-Beltrami operator $\Delta_g$.

Another famous equation similar to the prescribing curvature equation (\ref{pres-S2}) is the  mean field equation, namely
\be\label{meanfield}\Delta_g u=\rho\le(\f{he^u}{\int_\Sigma he^udv_g}-\f{1}{{\rm Vol}_g(\Sigma)}\ri),\ee
where $\rho$ is some real number. When $h$ is a smooth positive function, it was proved by Ding-Jost-Li-Wang \cite{DJLW} that for $\rho<8\pi$, (\ref{meanfield}) has a solution;
under certain geometric condition on $(\Sigma,g)$, (\ref{meanfield}) has a solution for $\rho=8\pi$ (when $\Sigma$ is a flat torus, this result was independently proved by
Nolasco-Tarantello \cite{NT}). The authors \cite{Yang-Zhu-Proc} derived the same conclusion when $h\geq 0$ and $h\not\equiv 0$.
Struwe-Tarantello \cite{Struwe} obtained a non-constant solution of (\ref{meanfield}), when $\Sigma$ is a
flat torus, $h\equiv1$ and $\rho\in(8\pi,4\pi^2)$.
When $\Sigma$ is a compact Riemannian surface of positive genus and $h$ is a smooth positive function,
Ding-Jost-Li-Wang \cite{DJLW2} proved that (\ref{meanfield}) admits a non-minimal solution for $\rho\in(8\pi,16\pi)$.
This result was generalized first by Chen-Lin \cite{ChenLin1, ChenLin2} to $\rho\in(8m\pi,16m\pi)$ ($m\in\mathbb{Z}^+$), then
by Malchiodi \cite{Malchiodi} to $\rho\in(8m\pi,16m\pi)$ ($m\in\mathbb{Z}^+$) and $\Sigma$ is a general Riemannian surface. Let $K$ be the Gaussian
curvature, Chen-Lin \cite{ChenLin1} also proved that, if $\D_g\log h(x)+\f{8m\pi}{{\rm Vol}_g(\Sigma)}-2K(x)>0$ for any $x\in\Sigma$,
then (\ref{meanfield}) has a solution. This improved the result in \cite{DJLW}.
By assuming $\Sigma$ to be a
flat torus and $h$ a
positive smooth function with certain symmetrization, Wang \cite{Wang} obtained an analog of \cite{DJLW}
for $\rho=16\pi$.

Our aim is to extend Ding-Jost-Li-Wang's result \cite{DJLW} to closed Riemannian surface with an isometric group action.
More precisely, let $(\Sigma,g)$, ${\textbf G}$ and $\ell$ be defined as in (\ref{Chen}). If $h$ is a smooth positive function
and $h(\sigma(x))=h(x)$ for all $\sigma\in
\textbf{G}$ and all $x\in\Sigma$, then for $\rho<8\pi\ell$, (\ref{meanfield}) has a solution; for $\rho=8\pi\ell$, under certain geometric hypothesis,
(\ref{meanfield}) has a solution.

\section{Notations and main results}

Let $(\Sigma,g)$ be a closed Riemannian surface and
$\textbf{G}=\{\sigma_1,\cdots,\sigma_N\}$ be
a finite isometric group acting on it. By definition, $\textbf{G}$ is a group and each $\sigma_i:\Sigma\ra\Sigma$
is an isometric map, particularly $\sigma_i^\ast g_x=g_{\sigma_i(x)}$ for all $x\in \Sigma$.
Let  $u:\Sigma\ra\mathbb{R}$ be a measurable function, we say that $u\in\mathscr{I}_{\textbf{G}}$ if $u$ is $\textbf{G}$-invariant,
namely $u(\sigma_i(x))=u(x)$ for any $1\leq i\leq N$ and almost every $x\in\Sigma$.
We denote $W^{1,2}(\Sigma,g)$
the closure of $C^\infty(\Sigma)$ under the norm
$$\|u\|_{W^{1,2}(\Sigma,g)}=\le(\int_\Sigma\le(|\nabla_gu|^2+u^2\ri)dv_g\ri)^{1/2},$$
where $\nabla_g$ and $dv_g$ stand for the gradient operator and the Riemannian volume element respectively.
Define a Hilbert space
\be\label{HG}\mathscr{H}_\textbf{G}= \le\{u\in W^{1,2}(\Sigma,g)\cap \mathscr{I}_\textbf{G}: \int_\Sigma udv_g=0\ri\}\ee
with an inner product
$$\langle u,v\rangle_{\mathscr{H}_{\textbf{G}}}=\int_\Sigma \langle\nabla_gu,\nabla_gv\rangle dv_g,$$
where $\langle\nabla_gu,\nabla_gv\rangle$ stands for the Riemannian inner product of $\nabla_gu$ and $\nabla_gv$.
Let $\Delta_g=-{\rm div}_g\nabla_g$ be the Laplace-Beltrami operator. For any $x\in\Sigma$, we set $I(x)=\sharp \textbf{G}(x)$, where $\sharp A$
stands for the number of all distinct points in the set $A$, and $\textbf{G}(x)=\{\sigma_1(x),\cdots,\sigma_N(x)\}$.
 Let
\be\label{ell}\ell=\inf_{x\in\Sigma}I(x).\ee
Clearly we have $1\leq\ell\leq N$ since $1\leq I(x)\leq N$ for all $x\in\Sigma$.
Chen's inequality (\ref{Chen}) is equivalent to the following: $\forall \beta\leq 4\pi\ell$, there holds
\be\label{T-M-ineq}\sup_{u\in\mathscr{H}_\textbf{G},\,\int_\Sigma|\nabla_gu|^2dv_g\leq 1}\int_\Sigma e^{\beta u^2}dv_g<\infty.\ee
Define a functional on $\mathscr{H}_\textbf{G}$ by
\be\label{functional}J_\gamma(u)=\f{1}{2}\int_\Sigma|\nabla_gu|^2dv_g-\gamma\log\int_\Sigma he^udv_g.\nonumber\ee
For any $0<\epsilon<1$, it follows from (\ref{T-M-ineq}) and a direct method of  variation that the minimizer $u_\epsilon$ of
the subcritical functional $J_{8\pi\ell(1-\epsilon)}$ exists. Denote $c_\epsilon=u_\epsilon(x_\epsilon)=\max_{\Sigma}u_\epsilon$.
As we shall see, if $c_\epsilon$ is bounded, then up to
a subsequence,
$u_\epsilon$ would converge to some $u_0$ as $\epsilon$ tends to zero, and $u_0$ is a minimizer of the critical functional $J_{8\pi\ell}$;
if $u_\epsilon$ blows up, i.e. $c_\epsilon\ra \infty$, then up to
a subsequence, $x_\epsilon$ converges to some $x_0$, $I(x_0)=\ell$, and  $u_\epsilon$ converges to certain Green function $\widetilde{G}_{x_0}$ satisfying
\be\label{Greenfunct}\le\{\begin{array}{lll}\Delta_g\widetilde{G}_{x_0}=8\pi\sum_{i=1}^{\ell}\delta_{\sigma_i(x_0)}
-\f{8\pi\ell}{{\rm Vol}_g(\Sigma)}\\[1.5ex] \int_\Sigma \widetilde{G}_{x_0}dv_g=0,\end{array}\ri.\nonumber\ee
where we assume without loss of generality $\textbf{G}(x_0)=\{\sigma_1(x_0),\cdots,\sigma_\ell(x_0)\}$.
In a normal coordinate system $\{x_1,x_2\}$ near $x_0$, the Green function $\widetilde{G}_{x_0}$ can be written as
\be\label{G-repres}\widetilde{G}_{x_0}(x)=-4\log r+\widetilde{A}_{x_0}+b_1x_1+b_2x_2+c_1x_1^2+2c_2x_1x_2+c_3x_2^2+O(r^3),\ee
where $\widetilde{A}_{x_0}$, $b_1$, $b_2$, $c_1$, $c_2$, $c_3$ are constants, $r(x)$ denotes the geodesic distance between $x$ and $x_0$.

Our first result reads

\begin{theorem}\label{Thm1} Let $(\Sigma,g)$ be a closed Riemannian surface, $\textbf{G}=\{\sigma_1,\cdots,\sigma_N\}$ be an isometric group acting on it,
and $\mathscr{H}_{\textbf{G}}$ be a function space defined as in (\ref{HG}). Let $\ell$ be defined as in (\ref{ell}).
Suppose $h$ is a smooth positive function on $\Sigma$
and $h(\sigma(x))=h(x)$ for all $\sigma\in\textbf{G}$ and all $x\in\Sigma$. Then
$(i)$ for any $\rho<8\pi\ell$, the equation (\ref{meanfield}) has a solution;
$(ii)$ if (\ref{meanfield}) has no solution for $\rho=8\pi\ell$, there holds
\be\label{lower}\inf_{u\in\mathscr{H}_{\textbf{G}}}J_{8\pi\ell}(u)\geq -4\pi\ell\max_{x\in\Sigma}
(2\log(\pi\ell h(x))+\widetilde{A}_x)-8\pi\ell,\ee
where $\widetilde{A}_x$ is defined as in  (\ref{G-repres}).
\end{theorem}

Noting that $J_{8\pi\ell}(0)\geq \inf_{u\in\mathscr{H}_{\textbf{G}}}J_{8\pi\ell}(u)$, we have the following:

\begin{theorem}\label{Thm2} Under the assumptions of Theorem \ref{Thm1}, if in addition
\be\label{cond}\log\int_\Sigma hdv_g>1+\f{1}{2}\max_{x\in\Sigma}
(2\log(\pi\ell h(x))+\widetilde{A}_x),\ee
 then (\ref{meanfield}) has a solution for $\rho=8\pi\ell$.
\end{theorem}

Later we shall give an example of $(\Sigma,g)$, $\textbf{G}$ and $h$ satisfying the hypothesis (\ref{cond}).
Assuming further $I(x)$ is a constant, we shall prove another existence result for the mean field equation (\ref{meanfield}),
namely

\begin{theorem}\label{Thm3}
Let $(\Sigma,g)$ be a closed Riemannian surface, $\textbf{G}=\{\sigma_1,\cdots,\sigma_\ell\}$ be an isometric group acting on it,
and $\mathscr{H}_{\textbf{G}}$ be a function space defined as in (\ref{HG}). Let $I(x)$ be defined as in (\ref{ell}).
Suppose that $I(x)\equiv \ell$ for all $x\in\Sigma$, that
$h$ is a smooth positive function satisfying $h(\sigma(x))=h(x)$ for all $\sigma\in \textbf{G}$ and all $x\in\Sigma$, that
$2\log(\pi\ell h(p))+\widetilde{A}_p=\max_{x\in\Sigma}
(2\log(\pi\ell h(x))+\widetilde{A}_x)$, and that in a normal coordinate system near $p$,
$$h(x)-h(p)=k_1x_1+k_2x_2+k_3x_1^2+2k_4x_1x_2+k_5x_2^2+O(r^3).$$
If in addition
\be\label{hy-2}\f{8\pi\ell}{{\rm Vol}_g(\Sigma)}-2K(p)+b_1^2+b_2^2-\f{\Delta_g h(p)}{h(p)}+\f{2(k_1b_1+k_2b_2)}{h(p)}>0,\ee
then (\ref{meanfield}) has a solution for $\rho=8\pi\ell$, where $K(p)$ denotes the Gaussian curvature of $(\Sigma,g)$ at $p$.
\end{theorem}

In the case $\textbf{G}=\{Id\}$, where $Id:\Sigma\ra \Sigma$ is the identity map, Theorems \ref{Thm1}-\ref{Thm3} are reduced to that of
Ding-Jost-Li-Wang \cite{DJLW}. Though we are in the spirit of Ding-Jost-Li-Wang \cite{DJLW} for the proof
of the general case of Theorems \ref{Thm1}-\ref{Thm3}, many technical difficulties need to be smoothed. The first issue is to construct
Green functions having many singular points by using elliptic estimates and the symmetric properties of $(\Sigma,g)$. The second one is to
derive a lower bound of $J_{8\pi\ell}$ on $\mathscr{H}_{\textbf{G}}$ by using the maximum principle. Here we use an argument of
our previous work \cite{Yang-Zhu-SCM}
which quite simplified that of \cite{DJLW}. The third issue is to construct test functions showing that (\ref{lower}) does not hold under
the hypothesis (\ref{hy-2}), which implies that (\ref{meanfield}) has a solution. For some technical reason we assume $I(x)$ is a constant
for all $x\in \Sigma$. Even in the case $\textbf{G}=\{Id\}$, our test function is different from that of Ding-Jost-Li-Wang \cite{DJLW}.\\

The remaining part of this paper is organized as follows: In Section \ref{Sec3}, we discuss existence and properties of
certain Green function in our setting; In Section \ref{Sec4}, we prove Theorems \ref{Thm1} and \ref{Thm2}; In Section \ref{Sec5},
we prove Theorem \ref{Thm3}; In Section \ref{Sec6}, we give examples of surfaces and isometric groups satisfying the hypotheses of
Theorems \ref{Thm2} and \ref{Thm3}. Throughout this paper, we do not distinguish sequence and subsequence.

\section{Green functions}\label{Sec3}
In this section, we discuss Green functions defined on a closed Riemannian surface $(\Sigma,g)$ with the action of an isometric group
$\textbf{G}=\{\sigma_1,\cdots,\sigma_N\}$.
Let $G(x,y)$ be the  usual Green function on $(\Sigma,g)$, which is a distributional solution to
\be\label{Gree}\Delta_{g,y}G(x,y)=8\pi\delta_x(y)-\f{8\pi}{{\rm Vol}_g(\Sigma)},\ee
where $\Delta_{g,y}$ is the Laplace-Beltrami operator with respect to the variable $y$, $\delta_x$ is the standard Dirac measure centered
at $x$.
It follows from (\cite{Aubin}, Theorem 4.13) that such a $G$ does exist and
can be normalized so that $\int_\Sigma G(x,y)dv_{g,y}=0$ for all $x\in \Sigma$. Moreover, $G(x,y)=G(y,x)$ and
$G\in C^\infty(\Sigma\times\Sigma\setminus\{(x,x): x\in\Sigma\})$. A key observation is the following:

\begin{proposition}\label{prop-1}
For all $\sigma\in\textbf{G}$ and all $x,y\in\Sigma$, $x\not=y$, there holds $G(\sigma(x),\sigma(y))=G(x,y)$.
\end{proposition}

\proof Take $\sigma\in\textbf{G}$. Replacing $x$ and $y$ by two points $\sigma(x)$ and $\sigma(y)$ in (\ref{Gree}), we have that
$$\Delta_{g,\,\sigma(y)}G(\sigma(x),\sigma(y))=8\pi\delta_{\sigma(x)}(\sigma(y))-\f{8\pi}{{\rm Vol}_g(\Sigma)}.$$
Since $\Delta_{g,\sigma(y)}=\Delta_{g,y}$ and $\delta_{\sigma(x)}(\sigma(y))=\delta_{x}(y)$, there holds  in the distributional sense
$$\Delta_{g,y}\le(G(\sigma(x),\sigma(y))-G(x,y)\ri)=0.$$
Hence $G(\sigma(x),\sigma(y))-G(x,y)\equiv C$ for some constant $C$. Since
$$\int_\Sigma G(\sigma(x),\sigma(y))dv_{g,\,\sigma(y)}=\int_\Sigma G(x,y)dv_{g,y}=0,$$
we conclude $C=0$ and the proposition. $\hfill\Box$\\

Given any $x\in\Sigma$. Assume that $I(x)=\sharp\textbf{G}(x)=
\sharp\{\sigma_1(x),\cdots,\sigma_N(x)\}=j$. Recall the definition (\ref{ell}),
$\ell=\inf_{x\in\Sigma}I(x)$. It follows that $\ell\leq j\leq N$. With no loss of generality, we can assume that
\be\label{tilde-j}\textbf{G}(x)=\{\sigma_1(x),\cdots,\sigma_j(x)\}.\ee
Then we have the following existence result:
\begin{proposition}\label{prop-2} Let $x\in \Sigma$ be fixed such that $I(x)=j$ and (\ref{tilde-j}) holds. Then
there exists a unique Green function $\widetilde{G}_x$ satisfying
$$\le\{\begin{array}{lll}\Delta_g\widetilde{G}_x=8\pi\sum_{i=1}^j\delta_{\sigma_i(x)}-\f{8\pi j}{{\rm Vol}_g(\Sigma)}\\[1.5ex]
\int_\Sigma \widetilde{G}_xdv_g=0.\end{array}\ri.$$
Moreover, $\widetilde{G}_x$ can be explicitly written as
\be\label{Gtilde}\widetilde{G}_x(y)=G(\sigma_1(x),y)+\cdots+G(\sigma_j(x),y),\nonumber\ee
where $G(x,y)$ be a solution of (\ref{Gree}) with $\int_\Sigma G(x,y)dv_{g,y}=0$.
\end{proposition}
\proof Note that in the distributional sense
$$\Delta_{g,y}\le(\widetilde{G}_x(y)-G(\sigma_1(x),y)-\cdots-G(\sigma_j(x),y)\ri)=0.$$
It follows that $\widetilde{G}_x(y)-G(\sigma_1(x),y)-\cdots-G(\sigma_j(x),y)\equiv C$ for some constant $C$.
Since $$\int_\Sigma \widetilde{G}_x(y) dv_{g,y}=\int_\Sigma {G}(\sigma_i(x),y)dv_{g,y}=0,\quad i=1,\cdots,j,$$
we have $C=0$, which is the desired result. $\hfill\Box$ \\

By elliptic estimates, one writes near $x$
\be\label{g-dec}G(x,y)=-4\log r+A_x+O(r),\ee
where $r$ denotes the geodesic distance between $x$ and $y$, $A_x=\lim_{r\ra 0}(G(x,y)+4\log r)$. In the same way, one can
write near $\sigma_i(x)$ \be\label{tild-g}\widetilde{G}_x=-4\log r+\widetilde{A}_{\sigma_i(x)}+O(r),\ee
where $r$ denotes the geodesic distance between $\sigma_i(x)$ and $y$, $\widetilde{A}_{\sigma_i(x)}=\lim_{r\ra 0}(\widetilde{G}_x(y)+4\log r)$.
Another key observation is the following:

\begin{proposition}\label{prop-3}
Let $x\in\Sigma$ be such that $I(x)=j$ and (\ref{tilde-j}) be satisfied. Assuming (\ref{g-dec}) and (\ref{tild-g}), we have  $\widetilde{A}_{\sigma_i(x)}=\widetilde{A}_{\sigma_1(x)}$, $i=1,\cdots,j$.
\end{proposition}

\proof By Proposition \ref{prop-2},
$$\widetilde{G}_x(y)=G(\sigma_1(x),y)+\cdots+G(\sigma_j(x),y).$$
Clearly
\be\label{s-1}
\widetilde{A}_{\sigma_i(x)}=A_{\sigma_i(x)}+\sum_{k\not=i,\,k=1}^jG(\sigma_i(x),\sigma_k(x)).
\ee
By Proposition \ref{prop-1}, we can see that
\bea\nonumber
A_{\sigma_i(x)}&=&\lim_{r\ra 0}(G(\sigma_i(x),y)+4\log r)\\\nonumber
&=&\lim_{r\ra 0}(G(\sigma_i(x),\sigma_i(z))+4\log r)\\\nonumber
&=&\lim_{r\ra 0}(G(x,z)+4\log r)\\\label{s-2}
&=&A_{\sigma_1(x)}.
\eea
Here in the second equality, $r$ denotes the geodesic distance between $\sigma_i(x)$ and $\sigma_i(z)$.
While in the third equality, $r$  denotes the geodesic distance between $x$ and $z$.
Since $G(\sigma_i(x),\sigma_k(x))=G(\sigma_k(x),\sigma_i(x))$, we conclude
\be\label{s-3}\sum_{k\not=i,\,k=1}^jG(\sigma_i(x),\sigma_k(x))=\sum_{k=2}^jG(\sigma_1(x),\sigma_k(x)),
\quad\forall i=1,\cdots,j.\ee
Inserting (\ref{s-2}) and (\ref{s-3}) into (\ref{s-1}), we have $\widetilde{A}_{\sigma_i(x)}=\widetilde{A}_{\sigma_1(x)}$,
$i=1,\cdots,j$, as desired. $\hfill\Box$

\section{Proof of Theorems \ref{Thm1} and \ref{Thm2} }\label{Sec4}

\subsection{Proof of Theorem \ref{Thm1}}
Clearly $(i)$ of Theorem \ref{Thm1} is implied by the following lemma.
\begin{lemma}\label{sub-critical}
For any $0<\epsilon<1$, there exists some $u_\epsilon\in\mathscr{H}_\textbf{G}$ such that
\be\label{sub}J_{8\pi\ell(1-\epsilon)}(u_\epsilon)=\inf_{u\in\mathscr{H}_\textbf{G}}J_{8\pi\ell(1-\epsilon)}(u).\ee
Moreover, $u_\epsilon$ is a solution of the equation
\be\label{E-L}\le\{\begin{array}{lll}
 \Delta_gu_\epsilon=8\pi\ell(1-\epsilon)\le(\lambda_\epsilon^{-1}he^{u_\epsilon}-\f{1}{{\rm Vol}_g(\Sigma)}\ri)\\[1.5ex]
 \lambda_\epsilon=\int_\Sigma he^{u_\epsilon}dv_g.
\end{array}\ri.\ee
\end{lemma}
\proof In view of (\ref{T-M-ineq}), a direct method of variation leads to the existence of $u_\epsilon\in \mathscr{H}_\textbf{G}$
satisfying (\ref{sub}). Then a straightforward calculation shows $u_\epsilon$ satisfies the Euler-Lagrange equation (\ref{E-L}). $\hfill\Box$\\

For the proof of $(ii)$ of Theorem \ref{Thm1}, we modify an argument of blow-up analysis performed by Ding-Jost-Li-Wang \cite{DJLW}.
Let $c_\epsilon=u_\epsilon(x_\epsilon)=\max_\Sigma u_\epsilon$. If $c_\epsilon$ is bounded, then applying elliptic estimates to
(\ref{E-L}), we conclude that up to a subsequence, $u_\epsilon\ra u_0$ in $C^1(\Sigma,g)$. It is easy to see that
$u_0$ is a solution of (\ref{meanfield}) for $\rho=8\pi\ell$. Therefore, under the assumption of $(ii)$ of Theorem \ref{Thm1}, there must hold
$c_\epsilon\ra +\infty$ as $\epsilon\ra 0$. With no loss of generality we assume
$x_\epsilon\ra x_0$ as $\epsilon\ra 0$.
By definition, we have for any $u\in\mathscr{H}_{\textbf{G}}$
$$J_{8\pi\ell(1-\epsilon)}(u)=\f{1}{2}\int_\Sigma|\nabla_gu|^2dv_g-8\pi\ell(1-\epsilon)\log
\int_\Sigma he^{u}dv_g.$$
It follows that for any $u\in\mathscr{H}_{\textbf{G}}$,
$J_{8\pi\ell}(u)=\lim_{\epsilon\ra 0}J_{8\pi\ell(1-\epsilon)}(u)\geq \lim_{\epsilon\ra 0}J_{8\pi\ell(1-\epsilon)}(u_\epsilon)$.
To prove $(ii)$ of Theorem \ref{Thm1}, it then suffices to prove
\be\label{low-bdd}\lim_{\epsilon\ra 0}J_{8\pi\ell(1-\epsilon)}(u_\epsilon)\geq -4\pi\ell\max_{x\in\Sigma}
(2\log(\pi\ell h(x))+\widetilde{A}_x)-8\pi\ell,\ee
where $\widetilde{A}_x$ is defined as in  (\ref{G-repres}).
Let \be\label{scal}r_\epsilon=\f{\sqrt{\lambda_\epsilon}}{\sqrt{8\pi\ell(1-\epsilon)h(x_0)}}e^{-c_\epsilon/2}\ee
and define a sequence of blow-up functions
\be\label{varphi}\varphi_\epsilon(y)={u}_\epsilon\le(\exp_{x_\epsilon}(r_\epsilon y)\ri)-c_\epsilon,\ee
where $\exp_{x_\epsilon}:T_{x_\epsilon}\Sigma\,(\cong\mathbb{R}^2)\ra \Sigma$ is an exponential map.
A straightforward calculation shows on the Euclidean ball $\mathbb{B}_{\delta r_\epsilon^{-1}}(0)$,
\bea\label{phi-equa}
\Delta_{g_\epsilon}\varphi_\epsilon(y)=\f{h(\exp_{x_\epsilon}(r_\epsilon y))}{h(x_0)}e^{\varphi_\epsilon(y)}
-\f{8\pi\ell(1-\epsilon)}{{\rm Vol}_g(\Sigma)}r_\epsilon^2.
\eea
Similar to (\cite{Yang-Zhu-SCM}, Lemma 2.7), we have the following:
\begin{lemma}\label{te-0}
Let $c_\epsilon$ and $r_\epsilon$ be as above. Then there exists some constant $C$ such that $r_\epsilon^2\leq Ce^{-\f{1}{2}c_\epsilon}$.
\end{lemma}
\proof Testing the equation (\ref{E-L}) by $u_\epsilon$, one gets
$$\|\nabla_gu_\epsilon\|_{L^2(\Sigma,g)}^2=\f{8\pi\ell(1-\epsilon)}{\lambda_\epsilon}\int_\Sigma hu_\epsilon
e^{u_\epsilon}dv_g\leq 8\pi\ell c_\epsilon.$$
This together with Chen's inequality (\ref{Chen}) leads to
$$\int_\Sigma he^{u_\epsilon}dv_g\leq C\int_\Sigma e^{4\pi\ell\f{u_\epsilon^2}{\|\nabla_gu_\epsilon\|_{L^2(\Sigma,g)}^2}+
\f{\|\nabla_gu_\epsilon\|_{L^2(\Sigma,g)}^2}{16\pi\ell}}dv_g\leq Ce^{
\f{\|\nabla_gu_\epsilon\|_{L^2(\Sigma,g)}^2}{16\pi\ell}}\leq Ce^{\f{1}{2}c_\epsilon}.$$
Hence we have by (\ref{scal}) that
$$r_\epsilon^2=\f{e^{-c_\epsilon}}{8\pi\ell(1-\epsilon)h(x_0)}\int_\Sigma he^{u_\epsilon}dv_g\leq Ce^{-\f{1}{2}c_\epsilon}.$$
This gives the desired result. $\hfill\Box$\\

The power of the above lemma is evident. It leads to $r_\epsilon\ra 0$ and a stronger estimate $r_\epsilon c_\epsilon^q\ra 0$ for
all $q>1$ as $\epsilon\ra 0$, which is very important during the process of blow-up analysis.
Then applying elliptic estimates to (\ref{E-L}), we obtain
\be\label{b-b}\varphi_\epsilon\ra \varphi\,\,\, {\rm in}\,\,\, C^1_{\rm loc}(\mathbb{R}^2)\,\,\,{\rm as}\,\,\, \epsilon\ra 0,\ee where
$\varphi$ is a distributional solution of
$$\le\{\begin{array}{lll}
-\Delta_{\mathbb{R}^2}\varphi(y)=e^{\varphi(y)}\quad {\rm in}\quad\mathbb{R}^2\\[1.5ex]
\int_{\mathbb{R}^2}e^{\varphi(y)}dy<\infty.\end{array}
\ri.$$
A classification theorem of Chen-Li \cite{C-L} gives
$\varphi(y)=-2\log(1+|y|^2/8)$ and thus
\be\label{bubble-energy}\int_{\mathbb{R}^2}e^{\varphi(y)}dy=8\pi.\ee

We now claim that
\be\label{concentration}\lim_{R\ra\infty}\lim_{\epsilon\ra 0}\lambda_\epsilon^{-1}\int_{B_{Rr_\epsilon}(x_\epsilon)}he^{u_\epsilon}dv_g=\f{1}{\ell}.\ee
In fact, by a change of variable, (\ref{scal}) and (\ref{b-b})
\bna
\lambda_\epsilon^{-1}\int_{B_{Rr_\epsilon}(x_\epsilon)}he^{u_\epsilon}dv_g&=&(1+o_\epsilon(1))\lambda_\epsilon^{-1}
\int_{\mathbb{B}_R(0)}h(\exp_{x_\epsilon}(r_\epsilon y))e^{u_\epsilon(\exp_{x_\epsilon}(r_\epsilon y))}r_\epsilon^2dy\\
&=&\le(\f{1}{8\pi\ell}+o_\epsilon(1)\ri)\int_{\mathbb{B}_R(0)}e^{\varphi_\epsilon(y)}dy\\
&=&\le(\f{1}{8\pi\ell}+o_\epsilon(1)\ri)\le(\int_{\mathbb{B}_R(0)}e^{\varphi(y)} dy+o_\epsilon(1)\ri).
\ena
This together with (\ref{bubble-energy}) leads to (\ref{concentration}).

Recalling $\ell=\inf_{x\in\Sigma}I(x)=\inf_{x\in\Sigma}\sharp\textbf{G}(x)$, we now calculate $I(x_0)=\sharp\textbf{G}(x_0)$ as below.
\begin{lemma}\label{con-ell}
$I(x_0)=\ell$.
\end{lemma}
\proof Denote $I(x_0)=k_0$. Suppose that $k_0>\ell$. With no loss of generality we assume $\textbf{G}(x_0)=\{\sigma_1(x_0),\cdots,\sigma_{k_0}(x_0)\}$.
Fix some $\delta>0$ sufficiently small such that
$$B_{\delta}(\sigma_i(x_0))\cap B_{\delta}(\sigma_j(x_0))=\varnothing,\,\,\,\forall\,1\leq i<j\leq k_0.$$
Note that $x_\epsilon\ra x_0$ and $r_\epsilon\ra 0$ as $\epsilon\ra 0$. For any fixed $R>0$, one gets $B_{Rr_\epsilon}(\sigma_i(x_\epsilon))\subset B_{\delta}(\sigma_i(x_0))$ for sufficiently small $\epsilon>0$.
Since
$h$ is $\textbf{G}$-invariant, we have by (\ref{concentration}) that
\be\label{con-i}\lim_{R\ra\infty}\lim_{\epsilon\ra 0}{\lambda_\epsilon^{-1}}\int_{B_{Rr_\epsilon}(\sigma_i(x_\epsilon))}he^{u_\epsilon}dv_g=\f{1}{\ell},
\,\,\,\forall\, i=1,\cdots,k_0.\ee
By (\ref{con-i}) and the assumption $k_0>\ell$, we obtain
$$1=\lambda_\epsilon^{-1}\int_\Sigma he^{u_\epsilon}dv_g\geq \sum_{i=1}^{k_0}\lambda_\epsilon^{-1}\int_{B_{Rr_\epsilon}(\sigma_i(x_\epsilon))}
he^{u_\epsilon}dv_g=\f{k_0}{\ell}>1,$$
which is impossible. Therefore $k_0=\ell$. $\hfill\Box$\\

Hereafter, by Lemma \ref{con-ell}, we assume $\sigma_i(x_0)\not=\sigma_j(x_0)$ for all $1\leq i<j\leq \ell$ and
$\textbf{G}(x_0)=\{\sigma_1(x_0),\cdots,\sigma_\ell(x_0)\}$. It follows from (\ref{con-i}) and the definition of $\lambda_\epsilon$ that
\be\label{t-0}\lim_{R\ra\infty}\lim_{\epsilon\ra 0}\lambda_\epsilon^{-1}\int_{\Sigma\setminus \cup_{i=1}^\ell
B_{Rr_\epsilon}(\sigma_i(x_\epsilon))}he^{u_\epsilon}dv_g=0.\ee
Since the functions on the righthand side of (\ref{E-L}) are bounded in $L^1(\Sigma,g)$, we have by (\cite{Yang-Zhu-SCM}, Lemma 2.10) that
$u_\epsilon$ is bounded in $W^{1,q}(\Sigma,g)$ for any $1<q<2$. As a consequence there exists some $\widetilde{G}_{x_0}\in W^{1,q}(\Sigma,g)$
such that for $0<r<2q/(2-q)$,
$u_\epsilon$ converges to $\widetilde{G}_{x_0}$ weakly in  $W^{1,q}(\Sigma)$, strongly in $L^r(\Sigma)$ and almost everywhere in $\Sigma$.
In view of (\ref{con-i}) and (\ref{t-0}), $\widetilde{G}_{x_0}$ satisfies
$$\Delta_g\widetilde{G}_{x_0}=8\pi\ell\le(\f{1}{\ell}\sum_{i=1}^{\ell}\delta_{\sigma_i(x_0)}-\f{1}{{\rm Vol}_g(\Sigma)}\ri)$$
in the distributional sense. Moreover, $\int_\Sigma \widetilde{G}_{x_0}dv_g=0$.
 Let $\widetilde{G}_{x_\epsilon}$ be a distributional solution of
 \begin{align*}
 \begin{cases}
      &\Delta_g\widetilde{G}_{x_\epsilon}=8\pi\le(\sum_{i=1}^{\ell}\delta_{\sigma_i(x_\epsilon)}-\f{\ell}{{\rm Vol}_g(\Sigma)}\ri) , \\
      & \int_\Sigma \widetilde{G}_{x_\epsilon}dv_g=0.
\end{cases}
\end{align*}
 The existence of  $\widetilde{G}_{x_\epsilon}$ is based on Proposition \ref{prop-2}.
 By elliptic estimates and Proposition \ref{prop-3}, $\widetilde{G}_{x_\epsilon}$ can be written as
 \be\label{G-rep}\widetilde{G}_{x_\epsilon}(x)=-4\log r+\widetilde{A}_{x_\epsilon}+O(r)\ee
 near $\sigma_i(x_\epsilon)$,  where $r$ denotes the geodesic distance between
 $x$ and $\sigma_i(x_\epsilon)$, $i=1,\cdots, \ell$. Similarly $\widetilde{G}_{x_0}$ can be represented by
 $\widetilde{G}_{x_0}(x)=-4\log r+\widetilde{A}_{x_0}+O(r)$ near $\sigma_i(x_0)$,  where $r$ denotes the geodesic distance between
 $x$ and $\sigma_i(x_0)$, $i=1,\cdots,\ell$. We now claim that
   \be\label{c-v}\widetilde{A}_{x_\epsilon}\ra \widetilde{A}_{x_0}\quad {\rm as}\quad \epsilon\ra 0.\ee
 To see this, with no loss of generality, we assume that $\sigma_1=Id$, the identity map. Write $A_{x_\epsilon}=\lim_{r\ra 0}
 (G_{x_\epsilon}(x)+4\log r)$ and $A_{x_0}=\lim_{r\ra 0}
 (G_{x_0}(x)+4\log r)$. By Proposition \ref{prop-2},
 $\widetilde{A}_{x_\epsilon}=A_{x_\epsilon}+\sum_{i=2}^\ell G(\sigma_i(x_\epsilon),x_\epsilon)$.
 Since $A_{x_\epsilon}\ra A_{x_0}$ and $x_\epsilon\ra x_0$ as $\epsilon\ra 0$, we immediately have
 $\widetilde{A}_{x_\epsilon}\ra \widetilde{A}_{x_0}=A_{x_0}+\sum_{i=2}^\ell G(\sigma_i(x_0),x_0)$ as
 $\epsilon\ra 0$. This confirms our claim (\ref{c-v}).

 One can easily see that
 $$\Delta_g(u_\epsilon-\widetilde{G}_{x_\epsilon})\geq 0\quad{\rm in}\quad \Sigma\setminus\cup_{i=1}^\ell B_{Rr_\epsilon}(\sigma_i(x_\epsilon)).$$
 In view of (\ref{b-b}), (\ref{G-rep}) and (\ref{c-v}), for all $1\leq i\leq \ell$, there holds
 $$(u_\epsilon-\widetilde{G}_{x_\epsilon})|_{\p B_{Rr_\epsilon}(x_\epsilon)}=-c_\epsilon+2\log\lambda_\epsilon-\widetilde{A}_{x_0}-
 2\log(\pi\ell h(x_0))+o_\epsilon(1)+o_R(1).$$
 By the maximum principle, we have on $\Sigma\setminus\cup_{i=1}^\ell B_{Rr_\epsilon}(\sigma_i(x_\epsilon))$
 \be\label{max-prin}u_\epsilon-\widetilde{G}_{x_\epsilon}\geq-c_\epsilon+2\log\lambda_\epsilon-\widetilde{A}_{x_0}-
 2\log(\pi\ell h(x_0))+o_\epsilon(1)+o_R(1).\ee
 By (\ref{b-b}), one calculates
 \bea\nonumber
 \int_{\cup_{i=1}^\ell B_{Rr_\epsilon}(\sigma_i(x_\epsilon))}|\nabla_gu_\epsilon|^2dv_g&=&\ell\int_{B_{Rr_\epsilon}(x_\epsilon)}
 |\nabla_gu_\epsilon|^2dv_g\\\nonumber
 &=&\ell\le(\int_{\mathbb{B}_R(0)}|\nabla \varphi|^2dx+o_\epsilon(1)\ri)\\\label{en-iner}
 &=&16\pi\ell\log\le(1+\f{R^2}{8}\ri)-16\pi\ell+o_\epsilon(1)+o_R(1),
 \eea
 where $o_\epsilon(1)\ra 0$ as $\epsilon\ra 0$ for any fixed $R>0$, $o_R(1)\ra 0$ as $R\ra \infty$.

 It follows from the divergence theorem and (\ref{E-L}) that
 \bea\nonumber
 \int_{\Sigma_{R,\epsilon}}|\nabla_gu_\epsilon|^2dv_g&=&
 \int_{\Sigma_{R,\epsilon}}u_\epsilon\Delta_gu_\epsilon dv_g+
 \int_{\p \Sigma_{R,\epsilon}}u_\epsilon\f{\p u_\epsilon}{\p \nu}ds_g\\\nonumber
 &=&8\pi\ell(1-\epsilon)\int_{\Sigma_{R,\epsilon}}\lambda_\epsilon^{-1}hu_\epsilon e^{u_\epsilon}dv_g
 -\f{8\pi\ell(1-\epsilon)}{{\rm Vol}_g(\Sigma)}\int_{\Sigma_{R,\epsilon}} u_\epsilon dv_g\\
 &&\quad+
 \int_{\p \Sigma_{R,\epsilon}}u_\epsilon\f{\p u_\epsilon}{\p \nu}ds_g,\label{a-0}
 \eea
  where  $\nu$ defined on the boundary $\p\Omega$ denotes the outer unit vector with respect to the domain $\Omega$ and we write for simplicity
 $$\Sigma_{R,\epsilon}=\Sigma\setminus\cup_{i=1}^\ell B_{Rr_\epsilon}(\sigma_i(x_\epsilon)).$$
 By (\ref{max-prin}),
 \bea
 &&8\pi\ell(1-\epsilon)\int_{\Sigma_{R,\epsilon}}\lambda_\epsilon^{-1}hu_\epsilon e^{u_\epsilon}dv_g\nonumber\\&&\geq
 8\pi\ell(1-\epsilon)\int_{\Sigma_{R,\epsilon}}\widetilde{G}_{x_\epsilon}\f{h e^{u_\epsilon}}{\lambda_\epsilon}dv_g
 -8\pi\ell(1-\epsilon)\int_{\Sigma_{R,\epsilon}}(c_\epsilon-2\log\lambda_\epsilon)\f{h e^{u_\epsilon}}{\lambda_\epsilon}dv_g\nonumber\\
 &&\quad-8\pi\ell(1-\epsilon)\int_{\Sigma_{R,\epsilon}}\le(\widetilde{A}_{x_0}+2\log(\pi\ell h(x_0))+o_\epsilon(1)+o_R(1)\ri)
 \f{h e^{u_\epsilon}}{\lambda_\epsilon}dv_g.\label{a-1}
 \eea
 By the divergence and (\ref{E-L}),
 \bea\nonumber
 8\pi\ell(1-\epsilon)\int_{\Sigma_{R,\epsilon}}\f{\widetilde{G}_{x_\epsilon}h e^{u_\epsilon}}{\lambda_\epsilon}dv_g
 &=&\int_{\Sigma_{R,\epsilon}}\widetilde{G}_{x_\epsilon}\le(\Delta_gu_\epsilon+\f{8\pi\ell(1-\epsilon)}{{\rm Vol}_g(\Sigma)}\ri)dv_g\\\nonumber
 &=&\int_{\Sigma_{R,\epsilon}}u_\epsilon\Delta_g\widetilde{G}_{x_\epsilon}dv_g-\int_{\p \Sigma_{R,\epsilon}}\widetilde{G}_{x_\epsilon}
 \f{\p u_\epsilon}{\p \nu}ds_g\\\nonumber
 &&+\int_{\p\Sigma_{R,\epsilon}}u_\epsilon\f{\p \widetilde{G}_{x_\epsilon}}{\p\nu}ds_g-\f{8\pi\ell(1-\epsilon)}{{\rm Vol}_g(\Sigma)}
 \int_{\cup_{i=1}^\ell B_{Rr_\epsilon}(\sigma_i(x_\epsilon))}\widetilde{G}_{x_\epsilon}dv_g\\\nonumber
 &=&-\f{8\pi\ell}{{\rm Vol}_g(\Sigma)}\int_{\Sigma_{R,\epsilon}}u_\epsilon dv_g-
 \int_{\p\Sigma_{R,\epsilon}}\widetilde{G}_{x_\epsilon}\f{\p u_\epsilon}{\p\nu}ds_g\\
 &&+\int_{\p \Sigma_{R,\epsilon}}u_\epsilon\f{\p \widetilde{G}_{x_\epsilon}}{\p\nu}ds_g
 -\f{8\pi\ell(1-\epsilon)}{{\rm Vol}_g(\Sigma)}\int_{\Sigma\setminus\Sigma_{R,\epsilon}}\widetilde{G}_{x_\epsilon}
 dv_g.
 \eea
 Also one has
 \bea
 \nonumber 8\pi\ell(1-\epsilon)\int_{\Sigma_{R,\epsilon}}(c_\epsilon-2\log\lambda_\epsilon)\f{he^{u_\epsilon}}{\lambda_\epsilon}dv_g&=&
 (c_\epsilon-2\log\lambda_\epsilon)\int_{\Sigma_{R,\epsilon}}\le(\Delta_gu_\epsilon+\f{8\pi\ell(1-\epsilon)}{{\rm Vol}_g(\Sigma)}\ri) dv_g\\
 \nonumber&=&(c_\epsilon-2\log\lambda_\epsilon)\le\{\int_{\cup_{i=1}^\ell\p B_{Rr_\epsilon(\sigma_i(x_\epsilon))}}\f{\p u_\epsilon}{\p \nu}ds_g\ri.\\
  &&\le.+\f{8\pi\ell(1-\epsilon)}{{\rm Vol}_g(\Sigma)}\le({\rm Vol}_g(\Sigma)-\ell{{\rm Vol}_g(B_{Rr_\epsilon}(x_\epsilon))}\ri)\ri\}
 \eea
 and in view of (\ref{t-0}),
 \be\label{o-1}
 8\pi\ell(1-\epsilon)\int_{\Sigma_{R,\epsilon}}(\widetilde{A}_{x_0}+2\log(\pi\ell h(x_0))+o_\epsilon(1)+o_R(1))\f{he^{u_\epsilon}}{\lambda_\epsilon}dv_g
 =o_\epsilon(1)+o_R(1).
 \ee
 Inserting (\ref{a-1})-(\ref{o-1}) into (\ref{a-0}), we have
 \bea\nonumber
 \int_{\Sigma_{R,\epsilon}}|\nabla_gu_\epsilon|^2dv_g&\geq&\f{8\pi\ell}{{\rm Vol}_g(\Sigma)}\int_{\cup_{i=1}^\ell B_{Rr_\epsilon}(x_\epsilon)}u_\epsilon dv_g
 -\f{8\pi\ell(1-\epsilon)}{{\rm Vol}_g(\Sigma)}\int_{\cup_{i=1}^\ell B_{Rr_\epsilon}(x_\epsilon)}\widetilde{G}_{x_\epsilon}dv_g\\\nonumber
 &&-\int_{\cup_{i=1}^\ell\p B_{Rr_\epsilon}(\sigma_i(x_\epsilon))}\f{\p u_\epsilon}{\p \nu}(u_\epsilon-\widetilde{G}_{x_\epsilon}+c_\epsilon-2\log\lambda_\epsilon)ds_g\\\nonumber
 &&-\int_{\cup_{i=1}^\ell \p B_{Rr_\epsilon}(\sigma_i(x_\epsilon))}u_\epsilon\f{\p\widetilde{G}_{x_\epsilon}}{\p\nu}ds_g+
 \f{8\pi\ell(1-\epsilon)}{{\rm Vol}_g(\Sigma)}\int_{\cup_{i=1}^\ell B_{Rr_\epsilon}(x_\epsilon)}u_\epsilon dv_g\\\nonumber
 &&+\f{8\pi\ell(1-\epsilon)}{{\rm Vol}_g(\Sigma)}(-c_\epsilon+2\log\lambda_\epsilon)\le({\rm Vol}_g(\Sigma)-\ell{{\rm Vol}_g(B_{Rr_\epsilon}(x_\epsilon))}\ri)\\\label{geq-1}
 &&+o_\epsilon(1)+o_R(1).
 \eea
 In view of (\ref{varphi}) and (\ref{b-b}), for $y\in\p \mathbb{B}_{R}(0)$, there holds
 \bna
 \f{\p u_\epsilon}{\p\nu}(\exp_{x_\epsilon}(r_\epsilon y))
 =\f{1}{r_\epsilon}\le(\f{\p \varphi}{\p r}(y)+o_\epsilon(1)\ri)
 =\f{1}{r_\epsilon}\le(-\f{R/2}{1+R^2/8}+o_\epsilon(1)\ri).
 \ena
 It follows that for $i=1,\cdots,\ell$
 $$\le.-\f{\p u_\epsilon}{\p \nu}\ri|_{\p B_{Rr\epsilon}(\sigma_i(x_\epsilon))}=\f{1}{r_\epsilon}\le(\f{R/2}{1+R^2/8}+o_\epsilon(1)\ri).$$
 This together with (\ref{max-prin}) leads to
 \bea\nonumber
 &&-\int_{\cup_{i=1}^\ell\p B_{Rr_\epsilon}(\sigma_i(x_\epsilon))}\f{\p u_\epsilon}{\p \nu}(u_\epsilon-\widetilde{G}_{x_\epsilon}+c_\epsilon-2\log\lambda_\epsilon)ds_g\\\label{b-0}
 &&\geq \f{\pi\ell R^2}{1+R^2/8}(-2\log(\pi\ell h(x_0))-\widetilde{A}_{x_0})+o_\epsilon(1)+o_R(1).
 \eea
 By (\ref{b-b}) and (\ref{G-rep}), we have
 \be\label{a-2}-\int_{\cup_{i=1}^\ell\p B_{Rr_\epsilon}(\sigma_i(x_\epsilon))}u_\epsilon\f{\p\widetilde{G}_{x_\epsilon}}{\p\nu}ds_g=
 8\pi\ell c_\epsilon-16\pi\ell \log\le(1+\f{R^2}{8}\ri)+o_\epsilon(1)+o_R(1).\ee
 It is easy to see that
 \be\label{a-3}\int_{\cup_{i=1}^\ell B_{Rr_\epsilon}(x_\epsilon)}\widetilde{G}_{x_\epsilon}dv_g=o_\epsilon(1),\quad
 \int_{\cup_{i=1}^\ell B_{Rr_\epsilon}(x_\epsilon)}u_\epsilon dv_g=o_\epsilon(1)\ee
 and \be\label{a-4}(-c_\epsilon+2\log\lambda_\epsilon){\rm Vol}_g(B_{Rr_\epsilon}(\sigma_i(x_\epsilon)))=o_\epsilon(1).\ee
 Inserting (\ref{b-0})-(\ref{a-4}) into (\ref{geq-1}), we have
 \bna
 \int_{\Sigma_{R,\epsilon}}|\nabla_gu_\epsilon|^2dv_g&\geq&8\pi\ell\le(-2\log(\pi\ell h(x_0))-\widetilde{A}_{x_0}\ri)
 +16\pi\ell(1-\epsilon)\log\lambda_\epsilon\\
 &&-16\pi\ell\log\le(1+\f{R^2}{8}\ri)+o_\epsilon(1)+o_R(1).
 \ena
 This together with (\ref{en-iner}) leads to
 \bea\nonumber
 \int_{\Sigma}|\nabla_gu_\epsilon|^2dv_g&\geq&8\pi\ell\le(-2\log(\pi\ell h(x_0))-\widetilde{A}_{x_0}\ri)
 +16\pi\ell(1-\epsilon)\log\lambda_\epsilon\\\label{g-e}
 &&-16\pi\ell+o_\epsilon(1)+o_R(1).
 \eea
 In view of (\ref{E-L}), it follows from (\ref{g-e}) that
 \bna
 J_{8\pi\ell(1-\epsilon)}(u_\epsilon)&=&\f{1}{2}\int_\Sigma|\nabla_gu_\epsilon|^2dv_g-8\pi\ell(1-\epsilon)\log\int_\Sigma he^{u_\epsilon}dv_g\\
 &\geq&4\pi\ell\le(-2\log(\pi\ell h(x_0))-\widetilde{A}_{x_0}\ri)
 -8\pi\ell+o_\epsilon(1)+o_R(1).
 \ena
 Therefore
 $$\lim_{\epsilon\ra 0}J_{8\pi\ell(1-\epsilon)}(u_\epsilon)\geq 4\pi\ell\le(-2\log(\pi\ell h(x_0))-\widetilde{A}_{x_0}\ri)
 -8\pi\ell,$$
 which immediately leads to (\ref{low-bdd}) and completes the proof of $(ii)$ of Theorem \ref{Thm1}. $\hfill\Box$

\subsection{Proof of Theorem \ref{Thm2}}
Suppose that (\ref{meanfield}) has no solution. By Theorem \ref{Thm1},
$$\inf_{u\in\mathscr{H}_{\textbf{G}}}J_{8\pi\ell}(u)\geq -4\pi\ell\max_{x\in\Sigma}
(2\log(\pi\ell h(x))+\widetilde{A}_x)-8\pi\ell,$$
where $\widetilde{A}_x$ is defined as in  (\ref{G-repres}). Hence
$$J_{8\pi\ell}(0)=-8\pi\ell\log\int_\Sigma hdv_g\geq -4\pi\ell\max_{x\in\Sigma}
(2\log(\pi\ell h(x))+\widetilde{A}_x)-8\pi\ell,$$
which contradicts (\ref{cond}). Therefore (\ref{meanfield}) has a solution. $\hfill\Box$

\section{Proof of Theorem \ref{Thm3}}\label{Sec5}
In this section, under the assumptions of Theorem \ref{Thm3}, we shall construct test functions to show
\be\label{test}\inf_{u\in\mathscr{H}_{\textbf{G}}}J_{8\pi\ell}(u)< -4\pi\ell\max_{x\in\Sigma}
(2\log(\pi\ell h(x))+\widetilde{A}_x)-8\pi\ell,\ee
where $\widetilde{A}_x$ is defined as in  (\ref{G-repres}). This together with Theorem \ref{Thm1} concludes
Theorem \ref{Thm3}. In the sequel,
we assume $\textbf{G}=\{\sigma_1,\cdots,\sigma_\ell\}$,
$I(x)=\sharp\textbf{G}(x)\equiv \ell$ for all $x\in\Sigma$.
Pick up some point $p\in \Sigma$ such that
\be\label{ma-x}\widetilde{A}_p+2\log(\pi\ell h(p))=\max_{x\in\Sigma}(\widetilde{A}_x+2\log(\pi\ell h(x))).\ee
It follows from $I(p)=\ell$ that  $\sigma_1(p),\cdots,\sigma_\ell(p)$ are different points on $\Sigma$.
Let \be\label{delta}\delta=\f{1}{4}\min\le\{{\rm inj}_g(\Sigma),\min_{1\leq i<j\leq \ell}d_g(\sigma_i(p),\sigma_j(p))\ri\}.\ee
Choose a normal coordinate system $(B_\delta(p),\exp_{p}^{-1};\{y^1,y^2\})$ near $p$.
More precisely, fixing an orthonormal basis $\{e_1,e_2\}$ of the tangent space $T_p\Sigma$, one can
write $y=(y^1,y^2)=\exp_p^{-1}(x)$ if the exponential map  $\exp_p: T_p\Sigma(\cong\mathbb{R}^2)\ra\Sigma$ maps
the tangent vector $y^1e_1+y^2e_2$ to
the point $x\in B_\delta(p)$.  One also denotes $\exp_p(y^1e_1+y^2e_2)$ by $\exp_p(y)$, and thus
$\exp_p^{-1}(B_\delta(p))=\mathbb{B}_\delta(0)$.
Let $\widetilde{G}_p$ be a Green function satisfying
\be\label{G-tilde}
\le\{\begin{array}{lll}
\Delta_g\widetilde{G}_p=8\pi\sum_{i=1}^\ell\delta_{\sigma_i(p)}-\f{8\pi\ell}{{\rm Vol}_g(\Sigma)}\\[1.5ex]
\int_\Sigma \widetilde{G}_pdv_g=0.
\end{array}\ri.
\ee
The existence of such a $\widetilde{G}_p$ is based on Proposition \ref{prop-1}. By Propositions \ref{prop-2} and
\ref{prop-3},  we have
 \be\label{G-inv}\widetilde{G}_p(\sigma_i(y))=\widetilde{G}_p(y),\quad i=1,\cdots,\ell.\ee
In the above mentioned normal coordinate system $(B_\delta(p),\exp_{p}^{-1};\{y^1,y^2\})$, by elliptic estimates,
 $\widetilde{G}_p$ can be written as
\be\label{Green-repres}
\widetilde{G}_p(\exp_p(y))=-4\log r+\widetilde{A}_p+\alpha(y)+\beta(y),
\ee
where $r=|y|=d_g(p,\exp_p(y))$, $\widetilde{A}_p$ is a constant, $\alpha(y)=b_1y^1+b_2y^2$, $\beta(y)\in C^1(\mathbb{B}_\delta(0))$ and
$$\beta(y)=c_1(y^1)^2+2c_2y^1y^2+c_3(y^2)^2+O(r^3).$$
To proceed, we have the following:
\begin{lemma}\label{c1-3-k} Let $c_1$ and $c_3$ are constants in (\ref{Green-repres}),
$K(p)$ be the Gaussian curvature of $(\Sigma,g)$ at $p$.  Then the following identity holds
$$c_1+c_3+\f{2}{3}K(p)=\f{4\pi\ell}{{\rm Vol}_g(\Sigma)}.$$
\end{lemma}

\proof We modify the argument of (\cite{DJLW}, Proposition 3.2). In a normal coordinate system near $p$,  the Riemannian metric
can be written as $g=dr^2+g^2(r,\theta)d\theta^2$ with
\be\label{g-rtheta}g(r,\theta)=r-\f{K(p)}{6}r^3+\phi(\cos\theta,\sin\theta)r^4+O(r^5),\ee
where $\phi(s,t)=\sum_{j=0}^3a_js^{3-j}t^j$ is a third order homogenous polynomial.

By the divergence theorem, we have for any $0<r<\delta$
\be\label{div}\int_{\cup_{i=1}^\ell\p B_r(\sigma_i(p))}\f{\p \widetilde{G}_p}{\p n}ds_g=\int_{\Sigma\setminus \cup_{i=1}^\ell B_r(\sigma_i(p))}
\Delta_g\widetilde{G}_p dv_g=-\f{8\pi\ell}{{\rm Vol}_g(\Sigma)}\int_{\Sigma\setminus \cup_{i=1}^\ell B_r(\sigma_i(p))}dv_g.\ee
Note that $\sigma^\ast g(x)=g(\sigma(x))$ for all $\sigma\in \textbf{G}$ and $x\in\Sigma$. In view of (\ref{G-inv}), (\ref{Green-repres}) and (\ref{g-rtheta}), one calculates
\bna
\int_{\cup_{i=1}^\ell\p B_r(\sigma_i(p))}\f{\p \widetilde{G}_p}{\p n}ds_g
&=&\ell\int_0^{2\pi}\le(-\f{4}{r}+2c_1r\cos^2\theta+4c_2r\cos\theta\sin\theta+2c_3r\sin^2\theta+O(r^2)\ri)\\
&&\quad\quad\quad\le(r-\f{K(p)}{6}r^3+O(r^4)\ri)d\theta\\
&=& -8\pi\ell+2\pi\ell\le(c_1+c_3+\f{2}{3}K(p)\ri)r^2+O(r^3)
\ena
and
$$-\f{8\pi\ell}{{\rm Vol}_g(\Sigma)}\int_{\Sigma\setminus \cup_{i=1}^\ell B_r(\sigma_i(p))}dv_g=-\f{8\pi\ell}{{\rm Vol}_g(\Sigma)}
\le({\rm Vol}_g(\Sigma)-\pi\ell r^2+O(r^4)\ri).$$
Inserting the above two estimates to (\ref{div}) and comparing the the terms involving $r^2$, we get the desired result. $\hfill\Box$\\

We define a sequence of functions $(\phi_\epsilon)_{\epsilon>0}$ by
$$\phi_\epsilon(x)=\le\{
\begin{array}{lll}
c-2\log\le(1+\f{r^2}{8\epsilon^2}\ri)+\widetilde{A}_p+\alpha(\exp_p^{-1}(\sigma_i^{-1}(x))),&x\in B_{R\epsilon}(\sigma_i(p)),\,\, i=1,\cdots,\ell\\[1.5ex]
\widetilde{G}_p(x)-\eta(\sigma_i^{-1}(x))\,\beta(\exp_p^{-1}(\sigma_i^{-1}(x))),&x\in B_{2R\epsilon}(\sigma_i(p))\setminus B_{R\epsilon}(\sigma_i(p))\\[1.5ex]
\widetilde{G}_p(x),&x\in\Sigma\setminus \cup_{i=1}^\ell B_{2R\epsilon}(\sigma_i(p)),
\end{array}\ri.$$
where $R\epsilon\ra 0$ as $\epsilon\ra 0$, $R$ and $c$ are constants depending only on $\epsilon$ and will be determined later, $r=r(x)$ denotes the geodesic distance between $x$ and $\sigma_i(p)$ for $x\in B_{R\epsilon}(\sigma_i(p))$,
$\eta\in C_0^\infty(B_{2R\epsilon}(p))$  satisfies $0\leq \eta\leq 1$, $\eta\equiv 1$ on $B_{R\epsilon}(p)$ and $|\nabla_g\eta|\leq 4/(R\epsilon)$,
both $\alpha$ and $\beta$ are functions defined as in (\ref{Green-repres}).

To ensure $\phi_\epsilon\in W^{1,2}(\Sigma,g)$, we require for all $x\in \p B_{R\epsilon}(\sigma_i(p))$, there holds
$$c-2\log\le(1+\f{r^2}{8\epsilon^2}\ri)+\widetilde{A}_p+\alpha(\exp_p^{-1}(\sigma_i^{-1}(x)))=\widetilde{G}_p(x)-\eta(\sigma_i(x)\,\beta(\exp_p^{-1}(\sigma_i^{-1}(x))).$$
This together with (\ref{Green-repres}) implies that
\be\label{c}c=2\log\le(1+\f{R^2}{8}\ri)-4\log (R\epsilon).\ee
One can easily check that
$$\phi_\epsilon(\sigma_i(x))=\phi_\epsilon(x),\quad\forall x\in\Sigma,\,\,\forall i=1,\cdots,\ell.$$
A straightforward calculation shows
\bea\nonumber
\int_{B_{R\epsilon}(\sigma_i(p))}|\nabla_g\phi_\epsilon|^2dv_g&=&\int_{B_{R\epsilon}(p)}|\nabla_g\phi_\epsilon|^2dv_g\\
\nonumber&=&\int_{B_{R\epsilon}(p)}\f{16r^2}{(r^2+8\epsilon^2)^2}dv_g+\int_{B_{R\epsilon}(p)}|\nabla_g(\alpha(\exp_p^{-1}(x)))|^2dv_g\\
&&\quad-
\int_{B_{R\epsilon}(p)}\f{8r\nabla_gr\nabla_g(\alpha(\exp_p^{-1}(x)))}{r^2+8\epsilon^2}dv_g.\label{energ-re}
\eea
Since $\int_0^{2\pi}\phi(\cos\theta,\sin\theta)d\theta=0$, we have
\bna
\int_{B_{R\epsilon}(p)}\f{16r^2}{(r^2+8\epsilon^2)^2}dv_g&=&\int_0^{2\pi}d\theta\int_0^{R\epsilon}\f{16r^2}{(r^2+8\epsilon^2)^2}
(r-\f{K(p)r^3}{6}+\phi(\cos\theta,\sin\theta)r^4+O(r^5))dr\\
&=&16\pi\log\le(1+\f{R^2}{8}\ri)-\f{16\pi R^2}{R^2+8}+\f{128}{3}\pi K(p)\epsilon^2\log\le(1+\f{R^2}{8}\ri)\\&&-\f{8}{3}\pi K(p)(R\epsilon)^2
+O(\epsilon^2)+O((R\epsilon)^4).
\ena
In a normal coordinate system near $p$, $g^{ij}(y)=\delta^{ij}+O(r^2)$, there holds
\bna
\int_{B_{R\epsilon}(p)}|\nabla_g(\alpha(\exp_p^{-1}(x)))|^2dv_g&=&\int_{\mathbb{B}_{R\epsilon}(0)}
g^{ij}(y)\p_i(b_1y^1+b_2y^2)\p_j(b_1y^1+b_2y^2)(1+O(r^2))dy\\
&=& \pi(b_1^2+b_2^2)(R\epsilon)^2+O((R\epsilon)^4).
\ena
Moreover,
\bna
\int_{B_{R\epsilon}(p)}\f{8r\nabla_gr\nabla_g(\alpha(\exp_p^{-1}(x)))}{r^2+8\epsilon^2}dv_g&=&4\int_0^{2\pi}d\theta
\int_0^{R\epsilon}\f{\p_rr^2\p_r(b_1r\cos\theta+b_2r\sin\theta)}{r^2+8\epsilon^2}g(r,\theta)dr\\
&=&O((R\epsilon)^4).
\ena
Inserting the above three estimates to (\ref{energ-re}), we have
\bea\nonumber
\int_{B_{R\epsilon}(\sigma_i(p))}|\nabla_g\phi_\epsilon|^2dv_g&=&16\pi\log\le(1+\f{R^2}{8}\ri)-\f{16\pi R^2}{R^2+8}+\f{128}{3}\pi K(p)\epsilon^2\log\le(1+\f{R^2}{8}\ri)\\&&-\f{8}{3}\pi K(p)(R\epsilon)^2+\pi(b_1^2+b_2^2)(R\epsilon)^2\nonumber\\&&
+O((R\epsilon)^4)+O(\epsilon^2)+O((R\epsilon)^4\log(1+{R^2}/{8})).\label{Repsilon}
\eea

Write $p_i=\sigma_i(p)$, $\eta_i=\eta\circ (\exp_p^{-1}\sigma_i^{-1})$, $\alpha_i=\alpha\circ \sigma_i^{-1}$ and $\beta_i=\beta\circ \sigma_i^{-1}$.
By the divergence theorem and (\ref{G-tilde}), we calculate
\bea\nonumber
\int_{\Sigma\setminus\cup_{i=1}^\ell B_{R\epsilon}(p_i)}|\nabla_g\phi_\epsilon|^2dv_g&=&
\int_{\Sigma\setminus\cup_{i=1}^\ell B_{R\epsilon}(p_i)}|\nabla_g\widetilde{G}_p|^2dv_g+
\int_{\cup_{i=1}^\ell B_{2R\epsilon}(p_i)\setminus B_{R\epsilon}(p_i)}|\nabla_g(\eta_i\beta_i)|^2dv_g\\\nonumber
&&-2\int_{\cup_{i=1}^\ell B_{2R\epsilon}(p_i)\setminus B_{R\epsilon}(p_i)}\nabla_g\widetilde{G}_p\nabla_g(\eta_i\beta_i)dv_g\\\nonumber
&=&-\int_{\cup_{i=1}^\ell\p B_{R\epsilon}(p_i)}\widetilde{G}_p\f{\p\widetilde{G}_p}{\p n}ds_g+
\int_{\Sigma\setminus\cup_{i=1}^\ell B_{R\epsilon}(p_i)}\widetilde{G}_p\Delta_g\widetilde{G}_p dv_g\\\nonumber
&&+\int_{\cup_{i=1}^\ell B_{2R\epsilon}(p_i)\setminus B_{R\epsilon}(p_i)}|\nabla_g(\eta_i\beta_i)|^2dv_g\\\label{out}
&&-2\int_{\cup_{i=1}^\ell B_{2R\epsilon}(p_i)\setminus B_{R\epsilon}(p_i)}\eta_i\beta_i\Delta_g\widetilde{G}_p dv_g+
2\int_{\cup_{i=1}^\ell\p B_{R\epsilon}(p_i)}\eta_i\beta_i\f{\p\widetilde{G}_p}{\p n}ds_g.\quad
\eea
By (\ref{Green-repres}),  we have in a normal polar coordinate system near $p$
\bea\nonumber \widetilde{G}_p(\exp_p(y))&=&-4\log r+\widetilde{A}_p+b_1r\cos\theta+b_2r\sin\theta+c_1r^2\cos^2\theta+2c_2r^2\cos\theta\sin\theta\\[1.2ex]
&&\quad+c_3r^2\sin^2\theta+\varrho(\cos\theta,\sin\theta)r^3+O(r^4),\label{G-polar}\eea
where $y=(y^1,y^2)=(r\cos\theta,r\sin\theta)$, $\varrho(s,t)$ is a third order homogenous polynomial with respect to $s$ and $t$. Obviously
\begin{align}\label{patial-G}
\f{\p\widetilde{G}_p}{\p n}=&-\f{4}{r}+b_1\cos\theta+b_2\sin\theta+2c_1r\cos^2\theta+4c_2r\cos\theta\sin\theta+2c_3r\sin^2\theta\nonumber\\
&+3r^2\varrho(\cos\theta,\sin\theta)+O(r^3).
\end{align}
Combining (\ref{g-rtheta}), (\ref{G-polar}) and (\ref{patial-G}), one has
\bea\nonumber
-\int_{\cup_{i=1}^\ell\p B_{R\epsilon}(p_i)}\widetilde{G}_p\f{\p\widetilde{G}_p}{\p n}ds_g&=&-\ell\int_
{\p B_{R\epsilon}(p)}\widetilde{G}_p\f{\p\widetilde{G}_p}{\p n}ds_g\\\nonumber
&=&-32\pi\ell\log(R\epsilon)+8\pi\ell \widetilde{A}_p+4\pi\ell(c_1+c_3)(R\epsilon)^2\\\nonumber
&&-\pi\ell(b_1^2+b_2^2)(R\epsilon)^2-2\pi\ell(c_1+c_3+\f{2}{3}K(p))\widetilde{A}_p(R\epsilon)^2\\\label{est-1}
&&+8\pi\ell(c_1+c_3+\f{2}{3}K(p))(R\epsilon)^2\log(R\epsilon)+O((R\epsilon)^4\log(R\epsilon)).
\eea
Similarly
\bea\nonumber
\int_{\Sigma\setminus\cup_{i=1}^\ell B_{R\epsilon}(p_i)}\widetilde{G}_p\Delta_g\widetilde{G}_p dv_g&=&
-\f{8\pi\ell}{{\rm Vol}_g(\Sigma)}\int_{\Sigma\setminus\cup_{i=1}^\ell B_{R\epsilon}(p_i)}\widetilde{G}_pdv_g
=\f{8\pi\ell}{{\rm Vol}_g(\Sigma)}\int_{\cup_{i=1}^\ell B_{R\epsilon}(p_i)}\widetilde{G}_pdv_g\\\nonumber
&=&-\f{32\pi^2\ell^2}{{\rm Vol}_g(\Sigma)}(R\epsilon)^2\log(R\epsilon)+\f{16\pi^2\ell^2}{{\rm Vol}_g(\Sigma)}(R\epsilon)^2
+\f{8\pi^2\ell^2}{{\rm Vol}_g(\Sigma)}\widetilde{A}_p(R\epsilon)^2\\&&\quad+O((R\epsilon)^4\log(R\epsilon)),\label{est-2}
\eea
\be\label{est-3}\int_{\cup_{i=1}^\ell B_{2R\epsilon}(p_i)\setminus B_{R\epsilon}(p_i)}|\nabla_g(\eta_i\beta_i)|^2dv_g=O((R\epsilon)^4),\ee
\be\label{est-4}
 -2\int_{\cup_{i=1}^\ell B_{2R\epsilon}(p_i)\setminus B_{R\epsilon}(p_i)}\eta_i\beta_i\Delta_g\widetilde{G}_p dv_g=
\f{16\pi\ell}{{\rm Vol}_g(\Sigma)}\int_{\cup_{i=1}^\ell B_{2R\epsilon}(p_i)\setminus B_{R\epsilon}(p_i)}\eta_i\beta_idv_g
=O((R\epsilon)^4)
\ee
and
\be\label{est-5}2\int_{\cup_{i=1}^\ell\p B_{R\epsilon}(p_i)}\eta_i\beta_i\f{\p\widetilde{G}_p}{\p n}ds_g=
-8\pi\ell(c_1+c_3)(R\epsilon)^2+O((R\epsilon)^4).\ee
Inserting (\ref{est-1})-(\ref{est-5}) to (\ref{out}), we obtain
\bna
\int_{\Sigma\setminus\cup_{i=1}^\ell B_{R\epsilon}(p_i)}|\nabla_g\phi_\epsilon|^2dv_g&=&-32\pi\ell\log(R\epsilon)+8\pi\ell \widetilde{A}_p+
4\pi\ell(c_1+c_3)(R\epsilon)^2\\
&&-\pi\ell(b_1^2+b_2^2)(R\epsilon)^2-2\pi\ell(c_1+c_3+\f{2}{3}K(p))\widetilde{A}_p(R\epsilon)^2\\&&+8\pi\ell(c_1+c_3+\f{2}{3}K(p))(R\epsilon)^2\log(R\epsilon)\\
&&-\f{32\pi^2\ell^2}{{\rm Vol}_g(\Sigma)}(R\epsilon)^2\log(R\epsilon)+\f{16\pi^2\ell^2}{{\rm Vol}_g(\Sigma)}(R\epsilon)^2
+\f{8\pi^2\ell^2}{{\rm Vol}_g(\Sigma)}\widetilde{A}_p(R\epsilon)^2\\
&&-8\pi\ell(c_1+c_3)(R\epsilon)^2+O((R\epsilon)^2\log(R\epsilon)).
\ena
This together with (\ref{Repsilon}) and Lemma \ref{c1-3-k} leads to
\bea\nonumber
\int_\Sigma|\nabla_g\phi_\epsilon|^2dv_g&=&16\pi\ell\log\le(1+\f{R^2}{8}\ri)-\f{16\pi\ell R^2}{R^2+8}+\f{128}{3}\pi\ell K(p)\epsilon^2
\log\le(1+\f{R^2}{8}\ri)\\&&-32\pi\ell\log(R\epsilon)+8\pi\ell \widetilde{A}_p+O((R\epsilon)^4\log(R\epsilon))+O(\epsilon^2)\nonumber\\
&&+O((R\epsilon)^4\log{1+\f{R^2}{8}}).\label{energy}
\eea

Now we estimate the average of the integral of $\phi_\epsilon$ as follows.
\bea\nonumber
\overline{\phi}_\epsilon&=&\f{1}{{\rm Vol}_g(\Sigma)}\int_\Sigma \phi_\epsilon dv_g\\\nonumber
&=&\f{1}{{\rm Vol}_g(\Sigma)}\le(\int_{\cup_{i=1}^\ell B_{R\epsilon}(p_i)}\le(c-2\log\le(1+\f{r^2}{8\epsilon^2}\ri)+\widetilde{A}_p+\alpha_i\ri)dv_g-
\int_{\cup_{i=1}^\ell B_{2R\epsilon}(p_i)\setminus B_{R\epsilon}(p_i)}\eta_i\beta_i dv_g\ri.\\\nonumber
&&\quad\le.+\int_{\Sigma\setminus\cup_{i=1}^\ell B_{2R\epsilon}(p_i)}\widetilde{G}_pdv_g+
\int_{\cup_{i=1}^\ell B_{2R\epsilon}(p_i)\setminus B_{R\epsilon}(p_i)}\widetilde{G}_pdv_g\ri)\\\label{average}
&=&-\f{16\pi\ell}{{\rm Vol}_g(\Sigma)}\epsilon^2\log\le(1+\f{R^2}{8}\ri)+O((R\epsilon)^4\log(R\epsilon))+O((R\epsilon)^4\log(1+R^2/8)).
\eea

Next we estimate $\int_\Sigma he^{\phi_\epsilon}dv_g$.
\bea\nonumber
\int_{B_{R\epsilon}(p_i)}e^{\phi_\epsilon}dv_g&=&\int_{B_{R\epsilon}(p)}e^{\phi_\epsilon}dv_g\\\nonumber
&=& e^{c+\widetilde{A}_p}\int_0^{2\pi}d\theta\int_0^{R\epsilon}\f{e^{b_1r\cos\theta+b_2r\sin\theta}}{(1+r^2/(8\epsilon^2))^2}
\le(r-\f{K(p)}{6}r^3+O(r^4)\ri)dr\\\nonumber
&=&\le(\f{1}{8}+\f{1}{R^2}\ri)^2e^{\widetilde{A}_p}\epsilon^{-2}\le\{8\pi \f{R^2}{R^2+8}-\f{32}{3}\pi K(p)\epsilon^2\log\le(1+\f{R^2}{8}\ri)\ri.\\
&&\qquad\qquad \le.+ 16\pi(b_1^2+b_2^2)\epsilon^2\log\le(1+\f{R^2}{8}\ri)+O(\epsilon^2)\ri\}.\label{integ-1}
\eea
Let $\delta$ be defined as in (\ref{delta}).
Note that for $y\in \mathbb{B}_\delta(0)$, there holds
\bna
e^{\widetilde{G}_p(\exp_p(y))}&=&r^{-4}e^{\widetilde{A}_p}\le\{1+b_1r\cos\theta+b_2r\sin\theta+c_1r^2\cos^2\theta+2c_2r^2\cos\theta\sin\theta+c_3r^2\sin^2\theta\ri.\\
&&\qquad\qquad+\f{1}{2}\le.(b_1^2r^2\cos^2\theta+2b_1b_2r^2\cos\theta\sin\theta+b_2^2r^2\sin^2\theta)+O(r^3)\ri\}.
\ena
By a straightforward calculation, we have
\bea\nonumber
\int_{B_\delta(p_i)\setminus B_{2R\epsilon}(p_i)}e^{\phi_\epsilon}dv_g&=&\int_{B_\delta(p)\setminus B_{2R\epsilon}(p)}
e^{\widetilde{G}_p}dv_g\\\nonumber
&=&\pi e^{\widetilde{A}_p}\le\{(2R\epsilon)^{-2}-(c_1+c_3+\f{b_1^2+b_2^2}{2})\log(R\epsilon)\ri.\\\label{integ-3}
&&\qquad\qquad+\le. \f{K(p)}{3}\log(R\epsilon)+O(1)\ri\}.
\eea
and
\bea\nonumber
\int_{B_{2R\epsilon}(p_i)\setminus B_{R\epsilon}(p_i)}e^{\phi_\epsilon}dv_g&=&
\int_{B_{2R\epsilon}(p)\setminus B_{R\epsilon}(p)}e^{\widetilde{G}_p-\eta_1\beta_1}dv_g\\
&=&\pi e^{\widetilde{A}_p}((R\epsilon)^{-2}-(2R\epsilon)^{-2}+O(1)).\label{integ-4}
\eea

For $y\in\mathbb{B}_{\delta}(0)$, we write
$$h(\exp_p(y))-h(p)=k_1r\cos\theta+k_2r\sin\theta+k_3r^2\cos^2\theta+2k_4r^2\cos\theta\sin\theta+k_5r^2\sin^2\theta+O(r^3).$$
Noting that $h(p_i)=h(p)$, $i=1,\cdots,\ell$, we have
\begin{align}\label{h-1}
&\int_{B_{R\epsilon}(p_i)}(h-h(p_i))e^{\phi_\epsilon}dv_g\nonumber\\
=&\int_{B_{R\epsilon}(p)}(h-h(p))e^{\phi_\epsilon}dv_g\nonumber\\
=&\le(\f{1}{8}+\f{1}{R^2}\ri)^2e^{\widetilde{A}_p}\epsilon^{-2}\le\{32\pi (k_3+k_5+k_1b_1+k_2b_2)\epsilon^2\log\le(1+{R^2}/{8}\ri)+O(\epsilon^2)\ri\}.
\end{align}
Similarly
\be\label{h-2}\int_{B_{2R\epsilon}(p_i)\setminus B_{R\epsilon}(p_i)}(h-h(p_i))e^{\phi_\epsilon}dv_g=
\int_{B_{2R\epsilon}(p)\setminus B_{R\epsilon}(p)}(h-h(p))e^{\phi_\epsilon}dv_g=O(1)\ee
and
\bea\nonumber\int_{B_{\delta}(p_i)\setminus B_{2R\epsilon}(p_i)}(h-h(p_i))e^{\phi_\epsilon}dv_g&=&
\int_{B_{\delta}(p)\setminus B_{2R\epsilon}(p)}(h-h(p))e^{\phi_\epsilon}dv_g\\\label{h-3}
&=&-\pi e^{\widetilde{A}_p}(k_3+k_5+k_1b_1+k_2b_2)\log{(R\epsilon)}+O(1).\eea
It follows from (\ref{integ-1}) and (\ref{h-1}) that
\bea\nonumber
\int_{B_{R\epsilon}(p_i)}he^{\phi_\epsilon}dv_g&=&h(p_i)\int_{B_{R\epsilon}(p_i)}e^{\phi_\epsilon}dv_g+
\int_{B_{R\epsilon}(p_i)}(h-h(p_i))e^{\phi_\epsilon}dv_g\\\nonumber
&=&\pi h(p)e^{\widetilde{A}_p}\epsilon^{-2}\le\{\f{R^2+8}{8R^2}
+16\le(b_1^2+b_2^2-\f{2}{3}K(p)+2\f{k_3+k_5+k_1b_1+k_2b_2}{h(p)}\ri)\ri.\\\label{ann-0}
&&\le.\times\le(\f{1}{8}+\f{1}{R^2}\ri)^2\epsilon^2\log\le(1+\f{R^2}{8}\ri)
+O(\epsilon^2)\ri\}.
\eea
In view of (\ref{integ-4}) and (\ref{h-2}),
\be\label{ann-1}\int_{B_{2R\epsilon}(p_i)\setminus B_{R\epsilon}(p_i)}he^{\phi_\epsilon}dv_g
=\pi h(p)e^{\widetilde{A}_p}((R\epsilon)^{-2}-(2R\epsilon)^{-2}+O(1)).\ee
Also we have by (\ref{integ-3}) and (\ref{h-3}) that
\bea\nonumber\int_{B_\delta(p_i)\setminus B_{2R\epsilon}(p_i)}he^{\phi_\epsilon}dv_g&=&
\pi h(p)e^{\widetilde{A}_p}\le\{(2R\epsilon)^{-2}-\le(c_1+c_3-\f{K(p)}{3}+\f{b_1^2+b_2^2}{2}\ri.\ri.\\
&&\quad\le.\le.+\f{k_3+k_5+k_1b_1+k_2b_2}{h(p)}\ri)\log(R\epsilon)+O(1)\ri\}.\label{ann-2}\eea
Clearly
$$\int_{\Sigma\setminus\cup_{i=1}^\ell B_\delta(p_i)}he^{\phi_\epsilon}dv_g=
\int_{\Sigma\setminus\cup_{i=1}^\ell B_\delta(p_i)}he^{\widetilde{G}_p}dv_g=O(1),$$
which together with (\ref{ann-0})-(\ref{ann-2}) leads to
\bea\nonumber
\int_\Sigma he^{\phi_\epsilon}dv_g&=&\f{\pi\ell h(p)e^{\widetilde{A}_p}}{8\epsilon^2}\le\{
1+\f{16}{R^2}+128\le(b_1^2+b_2^2-\f{2}{3}K(p)+2\f{k_3+k_5+k_1b_1+k_2b_2}{h(p)}\ri)\ri.\\\nonumber
&&\quad\quad\qquad\times\le(\f{1}{8}+\f{1}{R^2}\ri)^2\epsilon^2\log\le(1+\f{R^2}{8}\ri)-8\le(
c_1+c_3-\f{K(p)}{3}+\f{b_1^2+b_2^2}{2}\ri.\\\label{int-S}
&&\quad\quad\qquad \le.\le.+\f{k_3+k_5+k_1b_1+k_2b_2}{h(p)}\ri)\epsilon^2\log(R\epsilon)+O(\epsilon^2)\ri\}.
\eea
Now we choose $R$ such that $R^4\epsilon^2=\f{1}{\log(-\log{\epsilon})}$. It follows that $R^{-4}=o(\epsilon^2\log\epsilon)$
and $\log R/\log\epsilon=-1/2+o(1)$ as $\epsilon\ra 0$. This together with (\ref{int-S}) implies that
\bea\nonumber
\log\int_\Sigma he^{\phi_\epsilon}dv_g&=&\log\f{\pi\ell h(p)e^{\widetilde{A}_p}}{8\epsilon^2}+\f{16}{R^2}+4\le\{
b_1^2+b_2^2-\f{2}{3}K(p)+c_1+c_3+\ri.\\
&&\qquad\qquad \le.+2\f{k_3+k_5+k_1b_1+k_2b_2}{h(p)}+o_\epsilon(1)\ri\}\epsilon^2\log\f{1}{\epsilon}.\label{log-1}
\eea
In view of (\ref{energy}),
\bea\nonumber
\f{1}{2}\int_\Sigma|\nabla_g\phi_\epsilon|^2dv_g&=& 16\pi\ell\log\f{1}{\epsilon}-8\pi\ell\log 8-8\pi\ell+4\pi\ell \widetilde{A}_p+\f{128\pi\ell}{R^2}\\
&&\qquad\quad+\le(\f{64}{3}\pi\ell K(p)+o_\epsilon(1)\ri)\epsilon^2\log\f{1}{\epsilon}.
\label{energy-1}
\eea
Moreover, (\ref{average}) implies that
\be\label{average-1}\overline{\phi}_\epsilon=\le(-\f{16\pi\ell}{{\rm Vol}_g(\Sigma)}+o_\epsilon(1)\ri)\epsilon^2\log\f{1}{\epsilon}.\ee
Noting that $2(k_3+k_5)=-\Delta_gh(p)$ and combining (\ref{log-1})-(\ref{average-1}) and Lemma \ref{c1-3-k}, we conclude
\bea\nonumber
J_{8\pi\ell}(\phi_\epsilon-\overline{\phi}_\epsilon)
&=&\f{1}{2}\int_\Sigma|\nabla_g\phi_\epsilon|^2dv_g-8\pi\ell\log
\int_\Sigma he^{\phi_\epsilon}dv_g+8\pi\ell \overline{\phi}_\epsilon\\\nonumber
&=&-8\pi\ell-4\pi\ell \widetilde{A}_p-8\pi\ell \log(\pi\ell h(p))-32\pi\ell\le\{\f{8\pi\ell}{{\rm Vol}_g(\Sigma)}-2K(p)\ri.\\
&&\quad\le. +b_1^2+b_2^2-\f{\Delta_g h(p)}{h(p)}+\f{2(k_1b_1+k_2b_2)}{h(p)}+o_\epsilon(1)\ri\}\epsilon^2\log\f{1}{\epsilon}.\label{hy-c}
\eea
Under the hypothesis (\ref{hy-2}), one can see for (\ref{hy-c}) that
$$J_{8\pi\ell}(\phi_\epsilon-\overline{\phi}_\epsilon)< -4\pi\ell
(2\log(\pi\ell h(p))+\widetilde{A}_p)-8\pi\ell.$$
This together with (\ref{ma-x}) implies (\ref{test}) and completes the proof of Theorem \ref{Thm3}. $\hfill\Box$

\section{Geometric conditions}\label{Sec6}

In this section, we give examples of surfaces on which the assumptions of Theorems \ref{Thm2} and \ref{Thm3} may be satisfied.
As a consequence, on such surfaces the mean field equation (\ref{meanfield}) may have a solution for appropriate
$h$.

\subsection{Geometric condition for Theorem \ref{Thm2}}
For $x=(x_1,x_2)$, $y=(y_1,y_2)\in \mathbb{R}^2$,
$x\sim y$ if and only if $x_1=y_1+m$ and $x_2=y_2+k$ for some integers $m$ and $k$. Define $\tau:\mathbb{R}^2\ra \mathbb{T}^2=\mathbb{R}^2/\sim$
by $\tau(x)=[x]$, the equivalent class of $x$. If $U$ is an open subset of $\mathbb{R}^2$, then $\tau(U)$ is defined as an open subset
of the torus $\mathbb{T}^2$.
Let $g={(\tau^{-1})}^{\ast}(g_0)$ be a Riemannian metric on $\mathbb{T}^2$, where $g_0=d{x_1}^2+d{x_2}^2$ is the standard Euclidean metric. Then $(\mathbb{T}^2,g)$ is a flat torus.
Let $\tau_0: \mathbb{R}^2\ra\mathbb{R}^2$ be a transition given by $\tau_0(x_1,x_2)=(x_1+1/2,x_2)$ for all $(x_1,x_2)\in\mathbb{R}^2$,
We set $\textbf{G}=\{\sigma_1,\sigma_2\}$, where $\sigma_1(P)=P$ and
$\sigma_2(P)=\tau\circ\tau_0\circ\tau^{-1}(P)$ for all $P\in\mathbb{T}^2$. Then $\textbf{G}$ is an isometric group acting on $(\mathbb{T}^2,g)$.
Clearly $I(P)=\sharp \textbf{G}(P)=2$ for all $P\in\mathbb{T}^2$. Motivated by Ding-Jost-Li-Wang \cite{DJLW}, we define a function on $\mathbb{R}^2$ by
\bna\lambda(x)&=&4\pi\le({x_2}^2-x_2+\f{1}{6}\ri)-2\log\le(1+e^{-4\pi x_2}-2e^{-2\pi x_2}\cos(2\pi x_1)\ri)\\
&&-\sum_{n=1}^\infty\log
\le(1+e^{-4\pi n-4\pi  x_2}-2e^{-2\pi n-2\pi  x_2}\cos(2\pi x_1)\ri)\\&&-\sum_{n=1}^\infty\log
\le(1+e^{-4\pi n+4\pi  x_2}-2e^{-2\pi n+2\pi  x_2}\cos(2\pi x_1)\ri).\ena
According to (\cite{Lang}, pages 42-47), we have $\lambda(x_1+1,x_2)=\lambda(x_1,x_2)$ and
$\lambda(x_1,x_2+1)=\lambda(x_1,x_2)$. Moreover, for any fixed $P\in\mathbb{T}^2$, if we assume $0\in\tau^{-1}(P)$, then
$G(P,Q)=\lambda\circ \tau^{-1}(Q)$ is the usual Green function on $\mathbb{T}^2$. Precisely $G$ satisfies $G(P,Q)=G(Q,P)$ for all
$P,Q\in \mathbb{T}^2$ and
$$\le\{\begin{array}{lll}\Delta_{g,Q}G(P,Q)=8\pi\delta_P(Q)-8\pi\\[1.5ex]
\int_{\mathbb{T}^2}G(P,Q)dv_{g,Q}=0.\end{array}\ri.$$
It was proved by Ding-Jost-Li-Wang \cite{DJLW} that $G(P,Q)$ has the following decomposition near $P$
\be\label{totic}G(P,Q)=-4\log r+A_P+O(r^2),\nonumber \ee
where $r$ denotes the geodesic distance between $Q$ and $P$ and
\be\label{AP}A_P=-4\log (2\pi)+\f{2\pi}{3}-8\sum_{n=1}^\infty\log(1-e^{-2\pi n}).\nonumber\ee
Let $\widetilde{G}_P(Q)=G(P,Q)+G(\sigma_2(P),Q)$ for all $Q\in\mathbb{T}^2$. Then $\widetilde{G}_P$ is a $\textbf{G}$-invariant Green function
satisfying
$$\le\{\begin{array}{lll}\Delta_{g}\widetilde{G}_P(Q)=8\pi\le(\delta_P(Q)+\delta_{\sigma_2(P)}(Q)\ri)-16\pi\\[1.5ex]
\int_{\mathbb{T}^2}\widetilde{G}_P(Q)dv_{g}=0.\end{array}\ri.$$
By elliptic estimates, we can write near $P$
$$\widetilde{G}_P(Q)=-4\log r+\widetilde{A}_P+O(r).$$
Clearly
\be\label{tildeAP}\widetilde{A}_P=A_P+G(\sigma_2(P),P)=A_P+G(P,\sigma_2(P)).\ee

Using the inequality $\log(1+t)\geq t/(t+1)$ for $t>-1$, we have
\bea\nonumber
A_P+2+2\log\pi&\leq&2-4\log 2-2\log\pi+\f{2\pi}{3}+8\sum_{n=1}^\infty\f{e^{-2\pi n}}{1-e^{-2\pi n}}\\\nonumber
&\leq& 2-4\log 2-2\log\pi+\f{2\pi}{3}+\f{8}{1-e^{-2\pi}}\sum_{n=1}^\infty {e^{-2\pi n}}\\\nonumber
&=&2-4\log 2-2\log\pi+\f{2\pi}{3}+\f{8e^{-2\pi}}{(1-e^{-2\pi})^2}\\\label{approx-1}
&\approx& -0.9752.
\eea
Similarly
\bea\nonumber
G(P,\sigma_2(P))+2\log 2&=&\lambda(\f{1}{2},0)+2\log 2\\\nonumber
&=&\f{2\pi}{3}-2\log 2-8\sum_{n=1}^\infty\log(1+e^{-2\pi n})\\\nonumber
&\leq&\f{2\pi}{3}-2\log 2-8\sum_{n=1}^\infty \f{e^{-2\pi n}}{1+e^{-2\pi n}}\\\nonumber
&\leq&\f{2\pi}{3}-2\log 2-\f{8e^{-2\pi}}{1-e^{-4\pi}}\\\label{approx-2}
&\approx&0.6932.
\eea
In view of (\ref{tildeAP}), we have by combining (\ref{approx-1}) and (\ref{approx-2}) that
$\widetilde{A}_P<-2-2\log\pi-2\log 2$.
Since $\widetilde{A}_P$ is independent of the base point $P$,  we conclude
\be\label{maxim}\max_{P\in\mathbb{T}^2}\widetilde{A}_P<-2-2\log\pi-2\log 2.\ee
Since ${\rm Vol}_g(\mathbb{T}^2)=1$, (\ref{maxim}) is exactly (\ref{cond}) where $h$ is a constant and $\ell=2$.

In conclusion, we construct a flat torus $(\mathbb{T}^2,g)$ with the action of an isometric group $\textbf{G}=\{\sigma_1,\sigma_2\}$
such that (\ref{cond}) holds when $h$ is a positive constant. Set $h=c+\epsilon\phi$, where $c>0$ and $\epsilon>0$ are two constants,
$\phi:\mathbb{T}^2\ra \mathbb{R}$ is a smooth function satisfying $\phi(\sigma_1(x))=\phi(\sigma_2(x))$ for all $x\in\mathbb{T}^2$.
If $\epsilon$ is chosen sufficiently small, then (\ref{cond}) still holds.

\subsection{Geometric condition for Theorem \ref{Thm3}}

Let $\{\mathbf{a},\mathbf{b}\}$ be a basis on $\mathbb{R}^2$.
Consider a torus $\mathbb{T}^2=\mathbb{R}^2/\sim$, where $x\sim y$ if and only if $x-y=k\mathbf{a}$ for some integer $k$ or
$x-y=m\mathbf{b}$ for some integer $m$. Let $g$ be a flat metric induced by the standard Euclidean metric on $\mathbb{R}^2$.
Denote the Gaussian curvature of $(\mathbb{T}^2,g)$ by $K$. Then $K(x)\equiv 0$ for all $x\in\mathbb{T}^2$. For any fixed positive
integer $\ell$, we set $\sigma_1(x)=x$, $\sigma_2(x)=x+\f{1}{\ell}\mathbf{a}$, $\cdots$, $\sigma_\ell(x)=x+\f{\ell-1}{\ell}\mathbf{a}$
for all $x\in \mathbb{T}^2$. Obviously $\textbf{G}=\{\sigma_1,\cdots,\sigma_\ell\}$ is an isometric group acting on $(\mathbb{T}^2,g)$.
One can easily see that (\ref{hy-2}) holds when $h$ is a positive constant. As in the previous subsection, we set
$h=c+\epsilon\phi$, where $c>0$ and $\epsilon>0$ are two constants,
$\phi:\mathbb{T}^2\ra \mathbb{R}$ is a smooth function satisfying $\phi(\sigma(x))=\phi(x)$ for all $\sigma\in\textbf{G}$ and
all $x\in\mathbb{T}^2$.
If $\epsilon$ is chosen sufficiently small, then (\ref{hy-2}) still holds.\\

{\bf Acknowledgements}. This work is partly supported by the National Science Foundation of China (Grant Nos.
  11471014 and  11761131002).

\end{document}